# ROOT LAPLACIAN EIGENMAPS WITH THEIR APPLICATION IN SPECTRAL EMBEDDING

Shouvik Datta Choudhury

shouvikdc8645@gmail.com

shouvik@capsulelabs.in

**ABSTRACT**

The root laplacian operator or the square root of Laplacian which can be obtained in complete Riemannian manifolds in the Gromov sense has an analog in graph theory as a square root of graph-Laplacian. Some potential applications have been shown in geometric deep learning(spectral clustering) and graph signal processing.

## 1. INTRODUCTION

The heart of both our paper and Belkin and Niyogi's paper is the spectral theorem of operators like Laplacian which is a pseudodifferential operator of rank 0, the subset of which are integrodifferential operators. Pseudo-differential operators as well as integrodifferential operators are used in the microlocal analysis of algebraic geometry where they are used as a different formulation of D-modules and development of coherent sheaves. A microlocal analysis is a branch of mathematics that studies the behavior of solutions to partial differential equations (PDEs) near singular points, also known as characteristic points, in the domain where the PDE is defined. It is a powerful tool for understanding the properties of PDEs and the solutions they generate, and has applications in fields such as quantum mechanics, optics, and control theory. Techniques used in the microlocal analysis include Fourier integral operators, wavefront sets, and pseudodifferential operators. Sheaves, which are used as data structures in the theoretical formulation of distributed computing, topological signal processing, and topological data analysis – a branch of machine learning initiated by Gunnar Carlsson et.al., are loosely speaking a collection or a conglomeration of mapping or functions emanating from disjoint sets in topological manifolds after the construction of a "presheaf", so that the local data encoded in each sheaf can be transferred to global sheaves. In graph signal processing, the concept of sheaves and presheaves can be used to model the local and global variations of graph signals. Sheaves and presheaves are mathematical structures that allow us to study the behavior of functions defined on a topological space, such as a graph, in a way that is robust to changes in the underlying space.A presheaf on a graph is a collection of functions defined on the open sets of the graph, such as the neighborhoods of each vertex, that are consistent with the topology of the graph. The presheaf assigns to each vertex a function, called a local section, that describes the behavior of the signal in a neighborhood of that vertex. The presheaf also assigns to each open set of the graph a function, called a restriction, that describes how the local sections on different vertices are related.A sheaf is a presheaf that also satisfies the property of gluing, which states that the restrictions of the local sections on different vertices must agree on the intersections of their domains. This property ensures that the local sections on different vertices can be "glued" together to form a global function that is well-defined on the entire graph. Sheaves and presheaves provide a way to model the variations of a graph signal in a way that is robust to changes in the underlying graph. They allow us to study the behavior of the signal locally, at each vertex, and also globally, over the entire graph. This can be useful for tasks such as denoising and compression of graph signals, as well as clustering and classification of the vertices of the graph. In this paper, we strive for the uniqueness and existence of a certain type of operator, a "square root" of the Laplace operator which is provided by taking the positive square root of the eigenvalue of the Laplace operator, the existence of which is both provided by the commutative setting of Banach algebra and non-commutative setting of Von-Neumann algebra and its possible "location" in Hilbert space provided by the statement and proof of Kato's conjecture. Kato's square root conjecture is a mathematical conjecture that states that, for a large class of differential operators, the square root of the operator is a well-defined and bounded operator. The conjecture is named after the Japanese mathematician Tosio Kato, who first formulated it in the 1950s. The conjecture applies to differential operators that are symmetric and have a compact inverse. A differential operator is symmetric if it is equal to its adjoint, and it has a compact inverse if the inverse operator maps a dense subset of the domain to a compact subset of the range. Some examples of operators that satisfy these conditions are the Laplacian operator on a Riemannian manifold and the Schrödinger operator on a quantum mechanical system. The conjecture states that for these types of operators, the square root of the operator is a well-defined and bounded operator, meaning that it maps a dense subset of the domain to a dense subset of the range and its norm is finite. The existence, as well as the uniqueness of the square root, enables us to the possible implementation of the square root in the spectral embedding context of machine learning as well as the geometric deep learning scenario initiated by Bronstein et. al. The graph data structures form the backbone in the most optimistic scenario as well as the pragmatic scenario for geometric deep learning which by far attempts to encompass non-Euclidean data like Riemannian and Finslerian structures giving a far-reaching generalization of neural networks so that the different aspects of data or signals such as noise and outliers are understood much better than their Euclidean counterpart. The Dirichlet "energy" of a map or more specifically an embedding measures the Laplacian "energy" or the degree or the involvement of mapping described by the Laplacian or the Laplace-Beltrami operator for Riemannian manifolds in the continuous aspect and in the discrete spectrum the graph Laplacian which is obtained from the adjacency matrix " topology" of the graph. In the spectral graph theory scenario, the graph or the discrete Dirichlet energy measures the embedding effect of a "line" or a geodesic as every compact Riemannian manifold admits in all cases a closed geodesic as prescribed by the Hopf-Rinow theorem which is a fundamental theorem of manifolds irrespective any structure ( here "structure" designates metric spaces which has length measures in both Lesbeage and Riemann sense and also Borel measures.

A differential operator that is known as the Laplace operator is named after the influential mathematician Piere-Simon Laplace. A differential operator is an operator that has the backbone of a linear, non-linear, or even in the general case a multilinear mapping. This operator transforms or "morphs" a space both normed and un-normed in the context of various space theories such as Banach and Hilbert into another space under its relegated action, which is a proliferation of differentiation in simple terms. An integrodifferential is a combination of integration and differentiation, and one type of integrodifferential is a differential operator. Differential operators are a member of the larger class known as integrodifferential. They have a wide variety of applications in science and engineering, such as many mathematical and kernel transformations, such as the Laplace transform. They also have significant applications in differential geometry, namely in the study of well-known geometric flows such as the Ricci flow, the Yamabe flow, and its key derivatives, all of which have contributed to the study of the structure and behavior of manifolds. The diffeomorphic structure of mappings, which is inherent in manifolds, can also be enabled by Laplace operations. Manifolds are defined by their intrinsic structures. In addition, the Laplace operator plays an essential part in the development - of both an ensemble of string theory and a quantum field theory functional in physics. This role can be described as both crucial and pivotal. Laplace operator also forms a significant ensemble of many important building blocks of mathematics and physics, which is adequately demonstrated in the construction of continuous Laplacian flow, an emergent geometric flow. This flow is a good example of how Laplace operator plays an important role in these fields. The Laplace operator is used as an edge-detecting operator for an image when it is processed with image processing software. When operators are used to analyze an indeterminate or variable, it is frequently necessary to map it to another context or space so that it can be thoroughly analyzed for its properties, regardless of how subtle they may be. This is for the same reason that it is necessary to use a function or a mapping. As every physical occurrence and its inherent "churn" can be modeled by differential equations, this is the core idea of function, mapping, or operators. Among the majority of these equations, the Laplace operator represents an important feature of the majority of them. A "full-blown" extreme generalization is acquainted in topos theory and n-category theory, whose abstractedness brings a systematic and unified perspective. Extreme generalization is provided by category theory, in which every function or mapping in the usual context of the term is alternated by objects and their associated morphisms and other subtle constructs. The continuous Laplace operator also has its generalization in the Laplace Beltrami operator, which is valid for a manifold both in the Riemannian and Finslerian perspectives. This operator can be thought of as having two views: the Riemannian and the Finslerian. The category theoretical viewpoints of algebra are employed in a variety of situations, including those of theoretical computer science, Riemannian geometry, and various theoretical aspects of machine learning, in order to bring a deep framework for abstract analysis. Every physical phenomenon can be reduced to its noumenon, which is continuous underlying. The concept of continuity provides an effective and reliable support structure for differentiability in any space, as well as the subsequent expedition across the full mathematical spectrum. Developing acceptable representations for complex data is one of the core difficulties in machine learning and pattern recognition. Consideration is given to the issue of generating a representation for data existing on a low-dimensional manifold embedded in a high-dimensional space. We offer a geometrically driven approach for expressing high-dimensional data based on the correspondence between the root Laplacian, the Laplace Beltrami operator on the manifold, and the links to the heat equation. The algorithm presents a computationally efficient, locality-preserving method for nonlinear dimensionality reduction with a natural relationship to clustering. Several potential applications and examples are explored.

## 2. LITERATURE REVIEW AND SUBSEQUENT COMPARISON

In this study, we investigate a method for constructing a graph that incorporates neighborhood information from the data set. Using the concept of the root of graph Laplacian, we construct a low-dimensional representation of the data set that ideally maintains local neighborhood information. The representation map generated by the approach can be considered as a discrete approximation of a continuous map that comes naturally from the manifold's shape.

Several features of the algorithm and framework of the analysis described here are worth highlighting:

The fundamental algorithm is quite simple. Several local computations and one sparse eigenvalue problem are required. The solution represents the manifold's inherent geometric structure. However, a search for neighboring objects in a high-dimensional space is required. We see that there are a number of effective approximate methods for locating the nearest neighbors(e.g. Indyk, 2000). The algorithm is justified by the role of the Laplace Beltrami operator in delivering the best embedding for the manifold. The approximation of the manifold is provided by the adjacency graph derived from the data points. The Laplace Beltrami operator has connected approximately the correctly chosen weighted Laplacian of the adjacency graph. The central role of the Laplace Beltrami operator in the heat equation enables us to select the weight decay function in a systematic manner using the heat kernel. Therefore, the data embedding maps connect the eigenmaps of the Laplace Beltrami operator, which are intrinsically defined maps on the entire manifold. In order to comprehend dimensionality-reduction methods geometrically, the framework of this research makes explicit use of these linkages. Despite the fact that the connections between the Laplace Beltrami operator and the root-graph Laplacian are well-known to geometers and specialists in spectral graph theory (Chung, 1997; Chung, Grigor'yan, & Yau, 2000), we are not aware of any application to dimensionality reduction or data representation except Belkin and Niyogi's framework. Kondor and Lafferty (2002) were the first to apply diffusion kernels to graphs and other discrete structures. Because the root-Laplacian eigenmaps technique preserves locality, it is relatively insensitive to outliers and noise. Additionally, it is not susceptible to short-circuiting because only local lengths are utilised. We demonstrate that by attempting to maintain local information in the embedding, the method indirectly accentuates the data's natural clusters. The close ties to spectral clustering techniques discovered in learning and computer vision (namely, the methodology of Shi & Malik, 1997) are then readily apparent. This relationship between dimensionality reduction and clustering is examined in depth. In contrast, global approaches, such as those presented by Tenenbaum et al. (2000), seek to retain all pairwise geodesic distances between locations and do not exhibit any propensity to cluster. Since much of the discussion in Seung and Lee (2000), Roweis and Saul (2000), and Tenenbaum et al. (2000) is motivated by the potential role that nonlinear dimensionality reduction may play in human perception and learning, it is worthwhile to consider the implication of the previous remark in this context. From high-dimensional stimuli, the biological perceptual apparatus must extract low-dimensional structures. If the strategy of recovering such a low-dimensional structure is intrinsically local (as in the suggested algorithm), then a natural clustering will arise and may serve as the foundation for the creation of categories in biological perception. Since our method is based on the manifold's fundamental geometric structure, it is stable in terms of embedding. The representation will not change so long as the embedding is isometric. In the example with the moving camera, different camera resolutions (i.e., different selections of $n$ in the $n \times n$ picture grid) should result in embeddings of the same underlying manifold into spaces with vastly different dimensions. Our approach will generate equivalent representations regardless of the resolution.

The all-encompassing problem of dimensionality reduction is subsequently relegated. Provided a set $\mathbf{p}_1, \ldots, \mathbf{p}_k$ of $k$ points in $\mathbb{R}^l$, select a collection of points $\mathbf{q}_1, \ldots, \mathbf{q}_k$ in $\mathbb{R}^m (m \ll l)$ such that $\mathbf{q}_i$ "represents" $\mathbf{p}_i$. In this paper, we redirect to the designated case where $\mathbf{p}_1, \ldots, \mathbf{p}_k \in \mathcal{M}$ and $\mathcal{M}$ is a manifold embedded in $\mathbb{R}^l$. An optimal algorithm is now developed to controllably orchestrate a representative $\mathbf{q}_i$ 's for this designated case.

The possible difference between our framework and the aforesaid paper's framework is very subtle. The overall difference is that we are creating a new matrix that has been constructed using the square root of the eigenvector of the original matrix. It should have been mentioned beforehand all the square roots taken in this paper either of the matrix or of the eigenvectors are positive either by taking the positive square root of the positive operand. This new matrix has real entries but the root eigenvectors of the original matrix in the parlance of Singular Value Decomposition is the square root of the original matrix. It is unique because of the imposed positive-semidefiniteness of the original matrix which is graph Laplacian or the combinatorial Laplacian. This technique, if deployed, has many interesting applications in various aspects and in various inter-disciplinary scenarios which will be described later in this paper. Interestingly we are operating in the real mathematical field to avoid "duping" of data, which may result due to incorporation of excessive duplicate data if we work on a complex field which is hard to "dedupe".

## 3. ALGORITHM

Provided $k$ points $\mathbf{p}_1, \ldots, \mathbf{p}_k$ in $\mathbb{R}^l$, we contrive a weighted graph with $k$ nodes, one for each point, and a set of edges connecting neighboring points. The embedding map is now fabricated by computing the eigenvectors of the root of graph Laplacian. The algorithmic procedure is formally stated below.

1 Step 1 (Adjacency graph assembly). We contrive an edge between nodes $i$ and $j$ if $\mathbf{p}_i$ and $\mathbf{p}_j$ are in a close neighborhood. There are two perspectives: (a) $\epsilon -$ neighborhoods (parameter $\epsilon \in \mathbb{R}$ ). Nodes $i$ and $j$ are joined by an edge if $\|\mathbf{p}_i - \mathbf{p}_j\|^2 < \epsilon$ where the norm of the space is the Euclidean norm in $\mathbb{R}^l$.

(b) $n$ nearest neighbors (parameter $n \in \mathbb{N}$ ). Nodes $i$ and $j$ are joined by an edge if $i$ is among $n$ nearest neighbors of $j$ or $j$ is among $n$ nearest neighbors of $i$. This connection is pertinently symmetric.

2 Step 2 (Weights selection). [1] We have two archetypes for weighting the edges: (a) Heat kernel (parameter $t \in \mathbb{R}$ ). If nodes $i$ and $j$ are connected, put

$$W_{ij} = e^{-\frac{\|\mathbf{p}_i - \mathbf{p}_j\|^2}{t}}$$

otherwise, insert $W_{ij} = 0$. The rationalization for this selection of weights will be provided later. In a computer enactment of the algorithm, steps 1 and 2 are channelized simultaneously. (b) Trivial-case (no parameters ( $t$ is considered $\infty$ ) ). $W_{ij} = 1$ if vertices $i$ and $j$ are connected by an edge and $W_{ij} = 0$ if vertices $i$ and $j$ are not joined by an edge. This course of action discards the need to choose $t$.

3 Step 3 (eigenmaps construction). Surmising the graph $G$, contrived above, is connected. Otherwise, step 3 is continued for each connected component. Compute eigenvalues and eigenvectors for the generalized eigenvector problem,

$$\xi \mathbf{f} = \sqrt{|\lambda D \mathbf{f}|},$$

where $D$ is the diagonal weight matrix, and its entries are column (or row, since $W$ is symmetric) addition of $W$, $D_{ii} = \sum_j W_{ji}$. $L = D - W$ is the emblematic Laplacian matrix. Laplacian is a symmetric, positive semi-definite matrix that can be emblematic of an operator on functions defined on vertices of $G$. S, a positive semi-definite matrix which is a square root of Graph Laplacian matrix $L$. Its existence and uniqueness are aptly described in []. $\lambda$ have designated the eigenvalues of the Laplace operator. Considering $\mathbf{f}_0, \ldots, \mathbf{f}_{k-1}$ be the solutions of equation 2.1, ordered in accordance to their eigenvalues:

$$\xi \mathbf{f}_0 = \sqrt{|\lambda_0 D \mathbf{f}_0|}$$
$$\xi \mathbf{f}_1 = \sqrt{|\lambda_1 D \mathbf{f}_1|}$$
$$\xi \mathbf{f}_{k-1} = \sqrt{|\lambda_{k-1} D \mathbf{f}_{k-1}|}$$
$$0 = \sqrt{|\lambda_0|} \leq \sqrt{|\lambda_1|} \leq \cdots \leq \sqrt{|\lambda_{k-1}|}$$

We discard the eigenvector $\mathbf{f}_0$ corresponding to eigenvalue **o** and use the next $m$ eigenvectors for embedding in $m$-dimensional Euclidean space:

$$\mathbf{p}_i \rightarrow (\mathbf{f}_1(i), \ldots, \mathbf{f}_m(i))$$

## EMBEDDING QUALITY

The concept of dilation and distortion of mapping was introduced by the mathematician Mikhail Gromov in the context of metric spaces. These concepts are used to measure how much a mapping between two metric spaces stretches or distorts the distances between points. The dilation of a mapping between two metric spaces is a measure of how much the mapping expands the distances between points. It is defined as the worst-case ratio of the distance between the images of two points to the distance between the original points. More formally, if f is a mapping from a metric space (X,d) to a metric space (Y,d'), the dilation of f is defined as:

$$\text{dilation}(f) = \sup\{\frac{d'(f(x), f(y))}{d(x, y)}\}$$

where x and y are points in X, and d'(f(x),f(y)) is the distance between the images of x and y in Y. The dilation of a mapping is always greater than or equal to 1, with 1 indicating that the mapping preserves the distances between points and values greater than 1 indicating that the mapping expands the distances. The distortion of a mapping between two metric spaces is a measure of how much the mapping distorts the angles between vectors. It is defined as the worst-case ratio of the angle between the images of two vectors to the angle between the original vectors. The distortion of a mapping between two metric spaces is a measure of how much the mapping distorts the angles between vectors. It is defined as the worst-case ratio of the angle between the images of two vectors to the angle between the original vectors. More formally, if f is a mapping from a normed vector space (X,||.||) to a normed vector space (Y,||.||'), the distortion of f is defined as:

$$\text{distortion}(f) = \sup\{||f(v)||'/||v||/||f(w)||'/||w||\}$$

where v and w are vectors in X, and ||.||' is the norm in Y. The distortion of a mapping is always greater than or equal to 1, with 1 indicating that the mapping preserves the angles between vectors and values greater than 1 indicating that the mapping distorts the angles. Gromov used these concepts to study the large-scale geometry of metric spaces and to define the notion of Gromov-hyperbolic spaces. The dilation and distortion can be used to measure the large-scale properties of the mappings between metric spaces and to classify them according to their behavior. Gromov's work in this area has had a significant impact on the field of geometric group theory and the study of geometric properties of metric spaces.

The concept of distortion and the eigenvalues of a linear operator are closely related, particularly in the context of eigenvalue problems. Distortion is a measure of how much a linear operator distorts the angles between vectors, and eigenvalues are a measure of how much a linear operator scales the vectors. A linear operator is a function that maps a vector space to itself and preserves the linear structure of the space. The eigenvalues of a linear operator are scalars that satisfy the equation:

$$Av = \lambda v$$

where A is the operator, λ is the eigenvalue, and v is the eigenvector. The eigenvalues of a linear operator are used to study the properties of the operator, such as its stability, symmetry, and singularity. The distortion of a linear operator is related to the eigenvalues in the sense that it is determined by the ratio of the eigenvalues of the operator. More specifically, the distortion of a linear operator is the ratio of the largest eigenvalue to the smallest eigenvalue.

$$\text{Distortion}(A) = \lambda\_max/\lambda\_min$$

The quality of a signal processing operation is often measured by how well it preserves certain properties of the graph signal, such as the graph frequency content or the smoothness of the signal. The eigenvalues of the graph Laplacian matrix are closely related to the graph frequency content and can be used as a measure of the quality of a signal processing operation. One common quality parameter used in graph signal processing is spectral distortion, which measures how much a signal processing operation distorts the graph frequency content. The spectral distortion is defined as the ratio of the largest eigenvalue of the graph Laplacian matrix after the signal processing operation to the largest eigenvalue of the original graph Laplacian matrix.

Another quality parameter used in graph signal processing is the vertex-frequency distortion, which measures how much a signal processing operation distorts the smoothness of the signal. The vertex-frequency distortion is defined as the ratio of the largest absolute difference in the vertex frequency components of the processed signal to the largest absolute difference in the vertex frequency components of the original signal. There are other quality parameters used in graph signal processing such as energy distortion, which measures how much energy of the signal is preserved after a processing operation, and variance distortion which measures how much the variance of the signal is preserved after a processing operation. In summary, the eigenvalues of the graph Laplacian matrix are closely related to the graph frequency content and can be used as a measure of the quality of a signal processing operation in graph signal processing. The spectral distortion, vertex-frequency distortion, energy distortion, and variance distortion are some of the quality parameters used to measure how well a signal processing operation preserves certain properties of the graph signal such as the graph frequency content and the smoothness of the signal.

In graph signal processing, the signals are defined on the nodes of a graph rather than on a regular grid as in traditional signal processing. The goal of graph signal processing is to design algorithms that can process signals on graphs and extract useful information from them. There are several key operations that are commonly used in graph signal processing, such as filtering, denoising, compression, and classification. These operations are typically performed on the signals defined at the nodes of the graph.

It is clear that low eigenvalues keep in comparison to high value lower distortion if the values are positive real both in signal processing and graph embedding in manifold if the processing of signals at nodes is kept constant. A sparse matrix is a matrix that has a large number of zero elements. Sparse matrices have several advantages over dense matrices, which are matrices that have a relatively small number of zero elements. One advantage of sparse matrices is that they use less memory than dense matrices. Because sparse matrices have a large number of zero elements, they require less storage space than dense matrices with the same dimensions. This can be particularly beneficial when working with large matrices or when working with limited memory resources. Another advantage of sparse matrices is that they can be more efficient to operate on than dense matrices. Many numerical algorithms that involve matrices, such as linear algebra operations and matrix factorizations, can take advantage of the sparsity of the matrix to perform the calculations more efficiently. This can result in faster computation times and more efficient use of computational resources. Sparse matrices are also useful in graph signal processing, where the graph is represented by an adjacency matrix. In this context, sparse matrices are efficient to store and process, because they have fewer non-zero elements than dense matrices. Sparse matrices are also useful in machine learning and data mining, where the data is often high-dimensional and sparse, such as text data, images, and graphs. In this context, sparse matrices are used to represent the data and the models, and the sparsity of the matrix is exploited to reduce the complexity and the computation time of the models.

In summary, sparse matrices have several advantages over dense matrices. They use less memory, can be more efficient to operate on, and are useful in specific applications such as graph signal processing, machine learning, and data mining where the data is often high-dimensional and sparse.

# EXPLANATION OF THE ALGORITHM

The most common Sobolev space is the Sobolev space of order k, denoted by $\mathbf{H^k}$ which consists of functions that are k times weakly differentiable and whose k-th weak derivatives are square-integrable. The k-th weak derivative is a generalization of the classical k-th derivative, which allows us to include functions that are not necessarily smooth but have a certain degree of differentiability in a weak sense. Sobolev spaces are used in the study of PDEs as function spaces for solutions of equations. They also play a central role in the theory of elliptic and parabolic equations, as well as in the theory of Sobolev inequalities, which provide information on the regularity of solutions of differential equations. Sobolev spaces are also important in the study of partial differential equations on manifolds, where they provide a natural framework for the analysis of solutions to equations that are defined on curved spaces. They are also used in the study of geometric partial differential equations, such as the mean curvature flow and the minimal surface equation. Sobolev spaces form a natural and versatile tool to study partial differential equations and other areas of analysis, and they are widely used in mathematics, physics, engineering, and other fields.

4.1. Dirichlet energy. In mathematics, more specifically in partial differential equations, the Dirichlet energy gauges a function's variability. It is acquainted as a quadratic functional in the Sobolev space $H_1$, to be more detailed. The German mathematician Peter Gustav Lejeune Dirichlet is the originator of the Dirichlet energy, which is closely connected to Laplace's equation.

Definition 1. Provided an open set $\Omega \subseteq \mathbf{R}^n$ and a morphism $u: \Omega \to \mathbf{R}$ the Dirichlet energy of the function $\chi$ is the real number

$$E[\chi] = \frac{1}{2} \int_\Omega \| \nabla u(\chi) \|^2 \, dx$$

where $\nabla u: \Omega \to \mathbf{R}^n$ connotes the gradient vector field of the function $\chi$.

As it is the integration of a non-negative quantity, the Dirichlet energy is non-negative, i.e. $E[u] \geq 0$ for $u$. Solving Laplace's equation $-\Delta u(\chi) = 0$ for all $\chi \in \Omega$, subject to appropriate boundary conditions, is conforming to solving the variational problem of determining a function $u$ that satisfies the boundary conditions and has minimal Dirichlet energy. Such a solution is called a harmonic function and such solutions are the aspects of investigation in potential theory. In a more general scenario, where $\Omega \subseteq \mathbf{R}^n$ is replaced by any Riemannian manifold $M$, and $u: \Omega \to \mathbf{R}$ is replaced by $u: M \to \Phi$ for another (various) Riemannian manifold $\Phi$, the Dirichlet energy is given by the sigma model. The solutions to the Lagrange equations for the sigma model Lagrangian are those functions $u$ that minimize/maximize the Dirichlet energy. In this specific case of $u: \Omega \to \mathbf{R}$, it just shows that the Lagrange equations (or, equivalently, the HamiltonJacobi equations) provide the foundational requirements and accumulation for obtaining extremal solutions.

The minimization of the Dirichlet energy of continuous operators such as the Laplace-Beltrami operator and the graph Dirichlet energy in the discrete counterpart, which is a standard procedure for the minimization of graph constraints, is the foundation of the algorithm that was described above and which forms the basis of embedding and spectral clustering via Laplacian eigenmaps. These two processes are the basis of embedding and spectral clustering. In the instance of the square root operator, we are calculating the Dirichlet graph energy matrix in the discrete or in the sense of the graph matrix. This helps us reduce the energy further, which ultimately results in an embedding that is more stable than that achieved by Belkin and Niyogi. When taking the square root of the graph Laplacian matrix, however, there was a chance that complex values would arise. This problem might be solved by taking the modulus or absolute value of the respective function (Laplacian matrix ). In this environment, the existence of complex values does not have the appropriate significance that it should. Therefore, in order to bring clarity to this technical discourse, we may state that Laplaclain embedding is equivalent to the minimization of graph Dirichlet energy, and root Laplacian embedding is equivalent to the minimization of the square root of the graph Dirichlet energy or matrix. Both of these statements can be made in the context of the following: Laplaclain embedding is equivalent to the minimization of graph Dirichlet energy. The root will invariably refer to the square root of a matrix throughout the entirety of this paper. The existence of the square root of a positive semidefinite operator as well as its uniqueness for a positive semi-definite root as well as the existence of the square root of a positive semidefinite matrix as well as its uniqueness for a positive semi-definite root has both been rigorously obtained in (Horn, Johnson), as well as in (Nair, Singh). Equivalence can exist between matrices and operators in some circumstances, specifically those in which the mathematical field associated with the operator is real. In the following stages, the resultant methodology will have a deduction that is both noticeable and practicable. It's possible for a matrix to have multiple square roots, but a semi-definite Hermitian matrix will only ever have a single, distinct square root-even for the nth root.

To assist in the dimensionality-reduction procedure that is initiated by Belkin and Niyogi based on Laplacian eigenmaps, the datasets which were projected into target space are condensed.

Previously, the majority of manifolds considered were sub-manifolds of Euclidean space. This perpetuates the question of whether this is true in general, i.e., whether every smooth manifold can be embedded in some Euclidean space. H. Whitney responded positively to and only a weak form of Whitney's theorem is demonstrated in Theorem 2.

In geometry specifically, an embedding (or imbedding ) is one exemplification of some mathematical construct contained within another construct, such as a group that is a subgroup. When some object $X$ is said to be embedded in another object $Y$, the embedding is given by some injective and topology-preserving map $f: X \to Y$. The precise meaning of "topology-preserving" depends on the kind of mathematical structure of which $X$ and $Y$ are instances. In the connotations of category theory, a structure-preserving map is called a morphism. Given $X$ and $Y$, several different embeddings of $X$ in $Y$ may be possible. In many similar interesting cases, there is a standard (or "canonical") embedding, like those of the natural numbers in the integers, the integers in the rational numbers, the rational numbers in the real numbers, and the real numbers in the complex numbers., In such a scenario it is prevalent to detect the domain $X$ with its image $f(X)$ contained in $Y$, so that $f(X) \subseteq Y$.

Theorem 2. A compact smooth manifold can be imbedded in a Euclidean space.

The imbedding space's dimensions are outrageously large. For the projective space $P^m$, with the differentiable structure adjoined by $m$ charts this theorem generates an imbedding in $\mathbf{R}^{m(m+1)}$. In fact, Whitney's theorem asserts the following:

Theorem 3. Let $f: M \to N$ be a smooth map which is an imbedding on a closed subset $C \subset M$; let $\varepsilon$ be a continuous positive function on $M$. If $\dim M \geq 2\dim N + 1$, then there is an imbedding $g: M \to N$ $\varepsilon$-approximating $f$ and such that $f|C = g|C$.

## 5. OPTIMALITY

Let us first demonstrate that the embedding provided by the root-Laplacian eigenmaps algorithm desiccates and enhances local information optimality. Recall that given a data set, we construct a weighted graph $G = (V, E)$ with edges connecting nearby points to each other. For the purposes of this discussion, assume the graph is connected. Consider the problem of mapping the weighted graph $G$ to a line so that connected points stay as close together as possible. Let $\mathbf{q} = (q_1, q_2, \ldots, q_n)^T$ be such a map. A reasonable and foundational criterion for choosing a "quality" map is to reduce the following objective function,

$$\sqrt{\left| \sum_{ij} (q_i - q_j)^2 W_{ij} \right|}$$

under appropriate constraints. The objective function with our choice of weights $W_{ij}$ incurs a huge error if neighboring points $\mathbf{p}_i$ and $\mathbf{p}_j$ are mapped far apart. Therefore, minimizing it is an attempt to ensure that if $\mathbf{p}_i$ and $\mathbf{p}_j$ are "near," then $q_i$ and $q_j$ are close as well. The above-mentioned quantity under square-root sign is the graph Dirichlet energy. The retribution has been further reduced for taking the square root of Dirichlet energy in the discrete case in comparison to Laplace operator. It turns out that for any $q$, we have deduced from Belkin and Niyogi,

$$\sqrt{\left|\frac{1}{2}\sum_{i,j}(q_i-q_j)^2 W_{ij}\right|} = \sqrt{|\mathbf{q}^T L \mathbf{q}|}$$

where $L = D - W$ as a graph-theoretic definition. To see this, notice that $W_{ij}$ is symmetric and $D_{ii} = \sum_j W_{ij}$. Continuing the expression under square root is taken to be positive. Therefore, the minimization problem reduces to generating

$$\operatorname*{argmin}_y \sqrt{\mathbf{q}^T L \mathbf{q}}.$$

$$\frac{y}{\sqrt{\mathbf{q}^T D \mathbf{q}} = 1}$$

The constraint $\sqrt{\mathbf{q}^T D \mathbf{q}} = 1$ removes an "random" scaling factor in the embedding. Matrix $D$ provides a natural measure on the vertices of the graph. The bigger the value $D_{ii}$ (corresponding to the $i$ th vertex) is, the more "important" is that vertex. It pertains from equation 3.1 that $L$ is a positive semidefinite matrix, and the vector $y$ that minimizes the objective function is given by the minimum eigenvalue solution to the generalized eigenvalue problem:

$$\xi \mathbf{y} = \sqrt{\lambda D \mathbf{q}}$$

Here $\xi$ is the positive definite square root of $L$, the Laplacian matrix. Let 1 be the constant function taking 1 at each vertex. It is easy to see that **1** is an eigenvector with eigenvalue 0. If the graph is connected, **1** is the only eigenvector for $\lambda = 0$. To eliminate this trivial solution, which collapses all vertices of $G$ onto the real number 1, we put an additional constraint of orthogonality and look for

$$\operatorname{argmin} \sqrt{\mathbf{q}^T L \mathbf{q}}.$$
$$\sqrt{\mathbf{q}^T D \mathbf{q}} = 1$$
$$\sqrt{\mathbf{q}^T D \mathbf{1}} = 0$$

Thus, the solution is now given by the eigenvector with the least positive eigenvalue. The mathematical "state" $\mathbf{q}^T D \mathbf{1} = 0$ can be interpreted as removing a translation invariance in q. Now consider the more general problem of embedding the graph into $m^-$ dimensional Euclidean space. The embedding is given by the $k \times m$ matrix $Q = [\mathbf{q}_1 \mathbf{q}_2, \dots, \mathbf{q}_m]$, where the $i$ th row provides the embedding coordinates of the $i$ th vertex. Similarly, we need to minimize

$$\sqrt{\sum_{i,j} \|\mathbf{q}^{(i)} - \mathbf{q}^{(j)}\|^2 W_{ij}} = \sqrt{\operatorname{tr}(Q^T L Q)}$$

where $\mathbf{q}^{(i)} = [\mathbf{q}_1(i), \ldots, \mathbf{q}_m(i)]^T$ is the $m$-dimensional representation of the $i$ th vertex. This converges to finding

$$\operatorname*{argmin}_{Q^T D Q = I} \sqrt{\operatorname{tr}(Q^T L Q)}.$$

For the one-dimensional embedding problem, the constraint prevents collapse onto a point. For the $m$-dimensional embedding problem, the constraint presented above prevents collapse onto a subspace of dimension less than $m - 1$ ($m$ if, as in one dimensional case, we require orthogonality to the constant vector). Standard methods show that the solution is provided by the matrix of eigenvectors corresponding to the lowest eigenvalues of the generalized eigenvalue problem $\xi \mathbf{q} = \sqrt{\lambda D \mathbf{q}}$.

## 6. LAPLACE-BELTRAMI OPERATOR

The root-Laplacian of a graph is connected and coupled to the Laplace Beltrami operator on manifolds in a different perspective. Let $\mathcal{M}$ be a smooth, compact, $m$-dimensional Riemannian manifold. If the manifold is embedded in $\mathbb{R}^l$, the Riemannian edifice (metric tensor) on the manifold is spawn by the standard Riemannian edifice on $\mathbb{R}^l$. As we did with the graph, we are looking here for a map from the manifold to the real line such that points close together on the manifold are mapped close together on the line. Let $f$ be such a morphism. We can assume that $f: \mathcal{M} \to \mathbb{R}$ is doubly differentiable. Consider two neighboring points $\mathbf{p}, \mathbf{r} \in \mathcal{M}$. They are projected to $f(\mathbf{p})$ and $f(\mathbf{r})$ subsequently. Belkin and Niyogi demonstrated that

$$|f(\mathbf{p}) - f(\mathbf{r})| \leq \operatorname{dist}_{\mathcal{M}}(\mathbf{p}, \mathbf{r}) \, \|\nabla f(\mathbf{p})\| + o(\operatorname{dist}_{\mathcal{M}}(\mathbf{p}, \mathbf{r}))$$

The gradient $\nabla f(p)$ is a vector in the tangent space $T\mathcal{M}_p$, such that given another vector $\mathbf{v} \in T\mathcal{M}_p, df(\mathbf{v}) = \langle \nabla f(p), \mathbf{v} \rangle_{\mathcal{M}}$. Thus, we see that $\|\nabla f\|$ provides us with an estimate of how far apart $f$ maps nearby points. We, therefore, delve for a map that best preserves locality on approximation by trying to acquiesce

$$\operatorname*{argmin}_{\|f\|_{L^2(\mathcal{M})} = 1} \sqrt{\int_{\mathcal{M}} \|\nabla f(p)\|^2},$$

where the integral is taken with respect to the conforming measure on a Riemannian manifold. Note that minimizing $\sqrt{\int_{\mathcal{M}} \|\nabla f(x)\|^2}$ corresponds to minimizing $S\mathbf{f} = \sqrt{\frac{1}{2} \sum_{i,j} (f_i - f_j)^2 W_{ij}}$ on a graph. Here, $\mathbf{f}$ is a function on vertices, and $f_i$ is the value of $\mathbf{f}$ on the $i$ th node of the graph. It turns out that minimizing the objective function of the above-mentioned equation reduces to deduce eigenfunctions of the root of the Laplace Beltrami operator $\mathcal{L}$.

$$\sqrt{\int_{\mathcal{M}} \| \nabla f \|^2} = \sqrt{\int_{\mathcal{M}} \mathcal{L}(f)f} = \sqrt{\int_{\mathcal{M}} \xi^{\wedge}2(f)f} = \sqrt{\int_{\mathcal{M}} \lambda f^3}$$

Now,

$$\sqrt{\int_{\mathcal{M}} \lambda f^3} = \sqrt{\lambda} f^2 / \sqrt{2} = \int_{\mathcal{M}} \xi(f)f$$

We acclaim that $\xi$ is positive semidefinite. $f$ that minimizes $\sqrt{\int_{\mathcal{M}} \| \nabla f \|^2}$ has to be an eigen function of $\xi$. The spectrum of $\mathcal{L}$ on a compact manifold $\mathcal{M}$ is known to be discrete (Rosenberg, 1997). Hence the spectrum of $\mathcal{S}$ on a compact manifold $\mathcal{M}$ is also discrete. Let the eigenvalues (in increasing order) be $0 = \sqrt{\lambda}_0 \leq \sqrt{\lambda}_1 \leq \sqrt{\lambda}_2 \leq \cdots$, and let $f_i$ be the eigenfunction corresponding to eigenvalue $\lambda_i$. It is easily seen that $f_0$ is the constant function that maps the entire manifold to a single point. To avoid this eventuality, we require (just as in the graph setting) that the embedding map $f$ be orthogonal to $f_0$. It immediately deduces that $f_1$ is the optimal embedding map. Following the arguments of the previous section, we see that

$$\mathbf{p} \to (f_1(\mathbf{p}), \dots, f_m(\mathbf{p}))$$

provides the optimal $m$-dimensional embedding

## 7. DIAGRAMS OF SPECTRAL EMBEDDING:

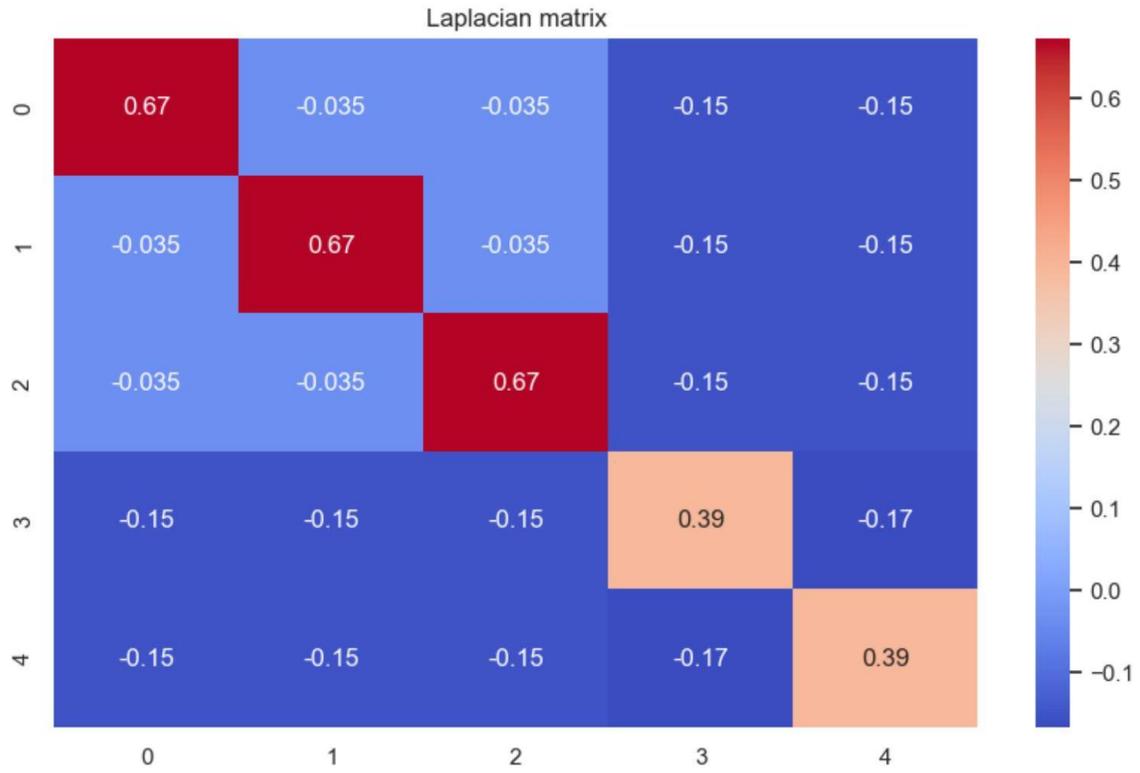

**Matrix heat map of Laplacian matrix of a random graph]**

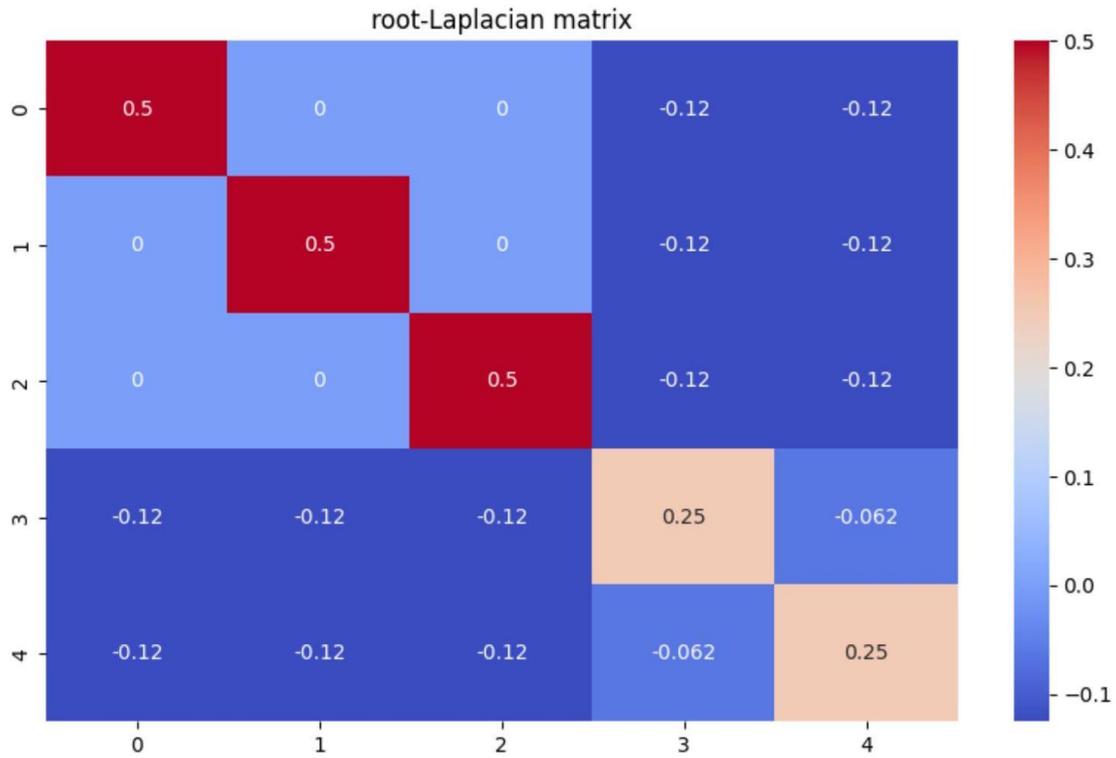

**Matrix heat map of root-Laplacian matrix of the same graph**

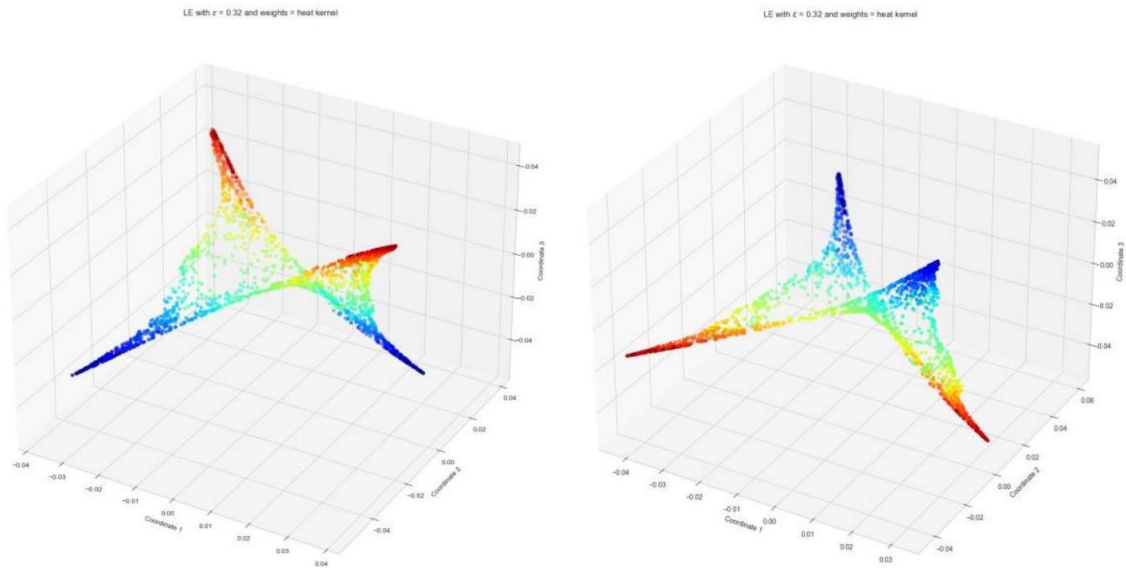

**10. [Spectral Embedding into 2D-dimension for the same heat kernel for laplacian and root laplacian respectively]**

# A POSSIBLE "IOT" Application for smart cities

The Internet of things is a collection or a network of "things" as the name suggests. Here "things" suggest are sensors, actuators, smart phones and other digital devices such as computational workstations both static and dynamic, an "improved" version of mobile computing treating all the components with some sort of "equality".

Smart cities are urban areas that use advanced technology and data analysis to improve the quality of life for citizens, increase the efficiency of city services, and reduce the environmental impact of the city.

A smart city typically uses a wide range of technologies such as the Internet of Things (IoT) sensors, big data analytics, cloud computing, and 5G networks to collect and analyze data from various sources such as traffic, energy consumption, weather, and air quality. The data is then used to optimize city services, improve decision-making and resource allocation, and enhance the overall livability and sustainability of the city.

Graph signal processing is a field of study that deals with the processing of signals defined on the vertices of a graph. The signals in this context are often referred to as graph signals. The graph structure can be used to represent the interactions between the different vertices, and the graph signal processing techniques can take advantage of this structure to extract meaningful information from the signals. Some common tasks in graph signal processing include filtering, denoising, and compression of graph signals, as well as clustering and classification of the vertices of the graph. Graph signal processing has applications in a wide range of fields including computer science, engineering, social sciences, and biology.

In graph signal processing, the Fourier basis is used to represent graph signals in the frequency domain. Just as the Fourier series is used to represent periodic functions as a sum of sine and cosine functions, the graph Fourier transform is used to represent graph signals as a sum of complex exponentials. Each vertex in the graph corresponds to a frequency, and the coefficients of the complex exponentials represent the amplitude of that frequency in the signal.

There are different ways to define the Fourier basis on a graph, but one common approach is to use the eigenvectors of the graph Laplacian matrix as the basis functions. The graph Laplacian is a matrix that encodes the structure of the graph and the eigenvectors of the Laplacian represent the different frequency modes of the graph signal.

The Graph Fourier Transform (GFT) is a linear operator that maps a signal defined on the vertices of a graph to its frequency domain representation. The GFT can be defined as the dot product of the signal with the eigenvectors of the Laplacian matrix.

Once the signal is represented in the frequency domain, it can be filtered, denoised, compressed, or manipulated. And it can be transformed back to the vertex domain by the Inverse Graph Fourier Transform(IGFT) which is the dot product of the signal with the eigenvectors of the Laplacian matrix.

This allows us to analyze and manipulate graph signals in a way that is similar to how we analyze and manipulate signals in traditional signal processing, by working in the frequency domain.

The most common Fourier basis is the set of complex exponentials, also known as sinusoids, which are defined as:

$$\int_{-\infty}^{+\infty} e^{2\pi i k t}$$

where t is the independent variable, k is an integer, and i is the imaginary unit. These functions are periodic with period 1 and form a complete orthonormal set, meaning that any function can be represented as a sum of these functions, and the sum converges to the original function in the mean square sense.

Another common Fourier basis is the set of trigonometric functions (sine and cosine) which are defined as:

$$\int_{-\infty}^{+\infty} \cos 2\pi i k t \sin 2\pi i k t$$

These functions are also periodic with period 1 and form a complete orthonormal set.

In graph signal processing, the Fourier basis is defined by the eigenvectors of the graph Laplacian matrix, which is a matrix that encodes the structure of the graph. The graph Laplacian is a symmetric matrix and its eigenvectors form an orthonormal basis for the space of functions defined on the graph.

The graph Laplacian matrix is defined as L = D - A, where D is the diagonal matrix of vertex degrees and A is the adjacency matrix of the graph. The eigenvectors of the Laplacian matrix are also called graph harmonics or graph modes.

Given a graph signal, which is a function defined on the vertices of the graph, the graph Fourier transform (GFT) is defined as the dot product of the signal with the eigenvectors of the Laplacian matrix. The dot product is used to calculate the projection of the signal on each eigenvector. The resulting coefficients are the frequency domain representation of the signal and are also called the graph Fourier coefficients.

The GFT and IGFT allow us to analyze and manipulate graph signals in a way that is similar to how we analyze and manipulate signals in traditional signal processing, by working in the frequency domain. In graph signal processing, the Fourier basis is defined by the eigenvectors of the graph Laplacian matrix, which is a matrix that encodes the structure of the graph. The graph Laplacian is a symmetric matrix and its eigenvectors form an orthonormal basis for the space of functions defined on the graph.

The graph Laplacian matrix is defined as L = D - A, where D is the diagonal matrix of vertex degrees and A is the adjacency matrix of the graph. The eigenvectors of the Laplacian matrix are also called graph harmonics or graph modes.

Given a graph signal, which is a function defined on the vertices of the graph, the graph Fourier transform (GFT) is defined as the dot product of the signal with the eigenvectors of the Laplacian matrix. The dot product is used to calculate the projection of the signal on each eigenvector. The resulting coefficients are the frequency domain representation of the signal and are also called the graph Fourier coefficients.

The concepts of frequency that come up in traditional signal processing offer a solid mathematical and logical foundation for signal analysis. Although it is technically possible, as was just discussed, to define notions of frequency for graph signals, it is not as easy to establish a corresponding understanding to comprehend these fundamental frequencies. It has been demonstrated empirically and supported conceptually that the frequency bases obtained from the shift operator tend to be ordered for the total variation criteria.

As of now, our attention has mainly been on frequency representations obtained from a graph's adjacency matrix. This strategy may be utilized with both directed and undirected graphs, and in the case of the cycle graph, it can be connected to DSP ideas. An undirected graph's Laplacian matrix can be used as the foundation for a frequency representation. Given that this matrix is positive semidefinite and has a complete set of orthogonal eigenvectors, all of the eigenvalues are real and non-negative, allowing us to write root GFT

$$\mathbf{S} = \sqrt[2]{|L|} = \sqrt{|\mathbf{U \Lambda U^\top}|}$$

$$\cong |\mathbf{U \Lambda^{\frac{1}{2}} U^\top}|$$

with U the GFT matrix of L, which is real-valued and orthogonal in this case. Because the eigenvalues are real, they provide a natural way to order the GFT basis vectors in terms of frequency (the variations of their values on the graph). In this case, the eigenvalue/eigenvector pairs can be viewed as successive optimizers of the Rayleigh quotient, where the $k$ th pair, $\lambda_k, \mathbf{u}_k$ solves :

$$\mathbf{u}_k = \sqrt{|\arg \min_{\mathbf{x}^\top \mathbf{u}_{k'}=0, k'=0,\ldots,k-1} \frac{\mathbf{x}^\top \mathbf{L} \mathbf{x}}{\mathbf{x}^\top \mathbf{x}}|}$$

with $\lambda_k = \sqrt{|\mathbf{u}_k^\top \mathbf{L} \mathbf{u}_{\mathbf{u}_k}|} \cong |\mathbf{U L^{\frac{1}{2}} U^\top}|$ which is equivalent if $\mathbf{u}_k$ is normalized. Thus for the explicit variation metric induced by the Laplacian quadratic form, the GFT provides an orthogonal basis with increased variation, such that, each additional basis vector minimizes the increase in variation while guaranteeing orthogonality. More generally, the relationships between eigenvectors and eigenvalues of the Laplacian and the structure of a graph are part of a deep and beautiful domain of mathematics known as spectral geometry.

The root-GFT and root-IGFT allow us to analyze and manipulate graph signals in a way that is similar to how we analyze and manipulate signals in traditional signal processing, by working in the frequency domain.

The combination of IoT and GSP can be used to improve the performance of various applications, such as:

Smart cities: IoT sensors can be used to collect data on various aspects of the city, such as traffic, energy consumption, and air quality. GSP can be used to analyze this data and extract useful information, such as traffic patterns, energy usage, and air quality trends.

Industrial IoT: IoT sensors can be used to collect data on industrial processes, such as temperature, humidity, and vibration. GSP can be used to analyze this data and extract useful information, such as equipment performance and maintenance needs.

Health monitoring: IoT devices can be used to collect data on various aspects of human health, such as heart rate, blood pressure, and sleep patterns. GSP can be used to analyze this data and extract useful information, such as patterns of health and wellness.

Overall, the combination of IoT and GSP can be used to improve the performance of various applications by collecting and analyzing data from IoT devices in a way that takes advantage of the graph structure of the data.